\documentclass[11pt,reqno]{amsart}

\usepackage{amsmath,amsfonts,amssymb}
\usepackage[latin1]{inputenc}
\usepackage{verbatim}
\usepackage{hyperref}
\usepackage{color}

\newtheorem{theorem}{Theorem}[section]
\newtheorem{lemma}[theorem]{Lemma}
\newtheorem{prop}[theorem]{Proposition}
\newtheorem{cor}[theorem]{Corollary}
\newtheorem{rmk}[theorem]{Remark}
\newtheorem{defn}[theorem]{Definition}

\newcommand{\sect}{\vspace{3mm} \setcounter{equation}{0} \section}

\newcommand{\w}[1]{\langle {#1} \rangle}

\newcommand{\pf}{\noindent {\bf Proof. \hspace{2mm}}}
\newcommand{\ef}{ \hfill $ \Box $ \vskip 3mm}
\newcommand{\be}{\begin{equation}}
\newcommand{\ee}{\end{equation}}
\newcommand{\bea}{\begin{eqnarray}}
\newcommand{\eea}{\end{eqnarray}}
\newcommand{\beas}{\begin{eqnarray*}}
\newcommand{\eeas}{\end{eqnarray*}}
\newcommand{\ep}{{\epsilon}}
\newcommand{\f}{\frac}
\newcommand{\na}{\nabla}

\newcommand{\bC}{{\mathbb C}}
\newcommand{\bR}{{\mathbb R}}
\newcommand{\bN}{{\mathbb N}}

\newcommand{\bS}{{\mathbb S}}

\newcommand{\vA}{{\mathcal A}}

\newcommand{\vD}{{\mathcal D}}

\newcommand{\vF}{{\mathcal F}}

\newcommand{\vH}{{\mathcal H}}
\newcommand{\vK}{{\mathcal K}}
\newcommand{\vL}{{\mathcal L}}

\newcommand{\vP}{{\mathcal P}}

\newcommand{\vT}{{\mathcal T}}

\newcommand{\vS}{{\mathcal S}}

\newcommand{\sgn}{\mbox{\rm sgn }}

\def\f{\frac}

\def\na{\nabla}

\def\ds{\displaystyle}

\def\blue{\color{blue}}

\begin{document}

\title[Resonances of the KFP operator]{
 Resonances of the Kramers-Fokker-Planck operator}

\author{Xue Ping WANG}
\address[X.P. Wang]{School of Mathematics and Key Laboratory of MIIT, Nanjing University of Aeronautics and Astronautics, Nanjing 211106, China, and Laboratoire de Math\'ematiques Jean Leray, Nantes Universit\'e, F-44322 Nantes Cedex 3, France}
\email{xue-ping.wang@univ-nantes.fr}


\begin{abstract}
We study the subellipticity of the analytically deformed Kramers-Fokker-Planck operator and give two definitions of resonances, one based on analytic dilation and the other on analytic distortion in the $x$-variables. Due to the failure of global subellipticity of the distorted operator in dimensions $n \geq 2$, the second method can only be used in dimension one. We show that, whenever both methods are applicable, they give rise to the same set of resonances. We also prove that, when $|\na V(x)| = O(\w x^{-\rho})$ with $\rho>1$, positive resonances are embedded eigenvalues and/or scattering resonances whose resonant states belong to weighted $L^2$-spaces and satisfy an outgoing radiation condition.
\end{abstract}

\maketitle


\sect{Introduction}
Resonances in quantum scattering are closely related to metastability and tunneling effects and they appear as poles of scattering matrices \cite{gam, wign}. The first mathematical definition of resonances for Schr\"odinger-type Hamiltonians in the scattering framework dates back to Aguilar and Combes \cite{ac}. They introduced the method of analytic dilation and considered resonances as discrete eigenvalues of a family of non-self-adjoint operators. Since then, the study of resonances has become an active research area at the intersection  of several branches of mathematics including partial differential equations, spectral theory and dynamical systems. See \cite{dz} for an overview. To date, most mathematical works on resonances have been concerned with self-adjoint elliptic operators.  In physics, metastability and tunneling effects also occur in the Kramers problem on the escape rate of a Brownian particle out of a deep potential well (\cite{kr,mel}). The associated Kramers-Fokker-Planck (KFP, for short)  operator is neither self-adjoint nor elliptic. The eigenvalues  in the low-temperature limit have been studied in \cite{hhs,hss} for external forces $F(x)$ satisfying a lower bound of the form $|F(x)| \ge C >0$ for $x$ outside some compact domain. The aim of this work is to establish a mathematical framework for resonances of the KFP operator when the force satisfies $|F(x)| \to 0$ as $|x|\to \infty$.
\\

The KFP equation, also called the Kramers equation by physicists (\cite{risk}),
 is the evolution equation for distribution functions describing  the Brownian motion of particles in an external field:
\be \label{kfp1}
\frac{\partial W}{\partial t} = \left( -v\cdot \nabla_x  + \nabla_v \cdot(\gamma v  -F(x))
+ \f{\gamma \kappa T }{m} \Delta_v \right)W,
\ee
where $W=W(t, x,v)$ for $x, v \in \bR^n$ and $t > 0$, and $F(x) = -  \nabla V(x)$ is the  external force.
Here  $x$ and $v$  represent  the position and velocity of particles, $m$ the mass, $\kappa$ the Boltzmann constant, $\gamma$ the friction coefficient and $T$ the temperature.  This equation is a special case of the more general  Fokker-Planck equation  or the Kolmogorov forward equation for continuous-time diffusion processes (\cite{risk}).

Set
\[
W(t, x, v) =e^{-\f{m}{2\kappa T} (\f{|v|^2}{2} + V(x))} u(t, x,v).
\]
Equation (\ref{kfp1}) is equivalent to
\be \label{kfp2}
\frac{\partial u}{\partial t} + \left( v\cdot \nabla_x  - \nabla V(x)\cdot \nabla_v  + \gamma (
- \f{\kappa T }{m} \Delta_v + \f{m }{4 \kappa T }|v|^2 -\f n 2) \right)u =0.
\ee
Let us omit the physical constants and denote
\be\label{operator} P=-\Delta_v+\f{1}{4}|v|^2-\f{n}{2} + v\cdot\na_x-\na V(x)\cdot\na_v, \quad (x,v) \in \bR^{2n}. \ee
$P$ is the KFP operator to be studied on $L^2(\bR^{2n}_{x,v})$. Throughout this work the potential
$V(x)$ is assumed to be at least a real-valued $C^1$ function on $\bR^n$.

Corresponding to the two different viewpoints, scattering and kinetic, we shall use the following two decompositions of $P$:
\be
P = P_0 + W= H_0 + L
\ee
where
\beas
P_0 &=& -\Delta_v+\f{1}{4}|v|^2-\f{n}{2} + v\cdot\na_x, \qquad W   =   -\na V(x)\cdot\na_v, \\
H_0  & =&   -\Delta_v+\f{1}{4}|v|^2-\f{n}{2}, \qquad
L  =   v\cdot\na_x-\na V(x)\cdot\na_v.
\eeas
By an abuse of notation, we still denote by $P$ ($P_0$, resp.) its closure with the maximal domain in  $L^2(\bR^{2n}_{x,v})$.\\

 The main  issue of this work is to establish global subelliptic estimates for the analytically deformed KFP operators, which  no longer satisfy H\"ormander's condition because the vector fields become complex when the deformation parameter $\theta$ is complex.  To use the method of analytic dilation, we assume that $V(x)$ is dilation-analytic, that is,
there exists  $a>0$ such that $V_\theta(x) = V(e^{\theta}x)$ defined for $\theta$ real admits a holomorphic extension in $\theta$ to the disc $\{\theta \in \bC, |\theta|<a\}$ and satisfies for some positive constants $\rho$ and $C$
\be \label{ass1}
|(\nabla_xV)(e^{\theta}x)| \le C \w{x}^{-\rho}, \qquad x\in \bR^n, |\theta| <a.
\ee
For example, $V(x) = \pm \w{x}^{r}$ for some $r <1$ satisfies  (\ref{ass1}) with $\rho=1-r>0$. Since we are interested in resonances with positive energies,   we can assume without loss of generality that $a\in ]0, \f \pi 4]$ (see Definition \ref{defn3.4} given below).
 For $\theta \in \bR$, let $U_\theta$ denote the unitary operator
\[
U_\theta f(x,v) = e^{\f{n\theta} 2} f(e^{\theta}x, v), \qquad f \in L^2(\bR^{2n}).
\]
Denote $P(\theta) = U_\theta P U_\theta^{-1}$, and similarly for $P_0(\theta)$, $W_\theta$ and $L_\theta$. Then
\be
P(\theta) = P_0(\theta) + W_\theta  = H_0 + L_\theta.
\ee
Under the assumption (\ref{ass1}), $P(\theta)$ defined for $\theta\in \bR$ can be extended to $\theta \in \bC$ with  $|\theta|<a$, and then to $|\Im \theta|<a$ by setting $P(\mu+i\varphi) = U_\mu P(i\varphi)U_\mu^{-1}$, $\mu\in\bR$,
$|\varphi|<a$; this is consistent, since $U_\mu P(\theta)U_\mu^{-1} = P(\theta+\mu)$ for $|\theta|<a$ and $|\theta+\mu|<a$ by
analytic continuation. All the results below for $|\Im\theta|<a$ follow from those for $|\theta|<a$ by this unitary
equivalence.

To define resonances of the KFP operator by analytic dilation, a key result is the following  global subelliptic estimate:
 for any $\phi_0\in ]0, \f \pi 2[$, there exist $ C>0$ and $s_0>0$ such that
\be \label{sub}
\| (H_0 +s)f\| +  \|\w{D_x}^{\f 2 3} f\|  + \| v\cdot \nabla_x f\| \le C \|({P}_0( \theta) +s) f\|,
\ee
 for $\theta \in \bC$ with $|\theta|\le \phi_0$,  $s \ge s_0$  and $f \in \vS(\bR_{x,v}^{2n})$  (see Theorem \ref{th2.1}). Estimate (\ref{sub}) is well known for $\theta=0$ (see, for example, \cite{n}).
 Its proof for $\Im \theta \neq 0$ is given in Section 2, by means of an explicit formula for the integral kernel of the semigroup generated by the fibre operators of $-e^{i\Im\theta}P_0(\theta)$ (see \cite{wz} for $\theta =0$). This method can also be used to prove $L^p$-subelliptic estimates (Theorem \ref{th-Lp}).  Estimate (\ref{sub}) shows that under the condition  (\ref{ass1}), $W_\theta$ is a relatively compact perturbation of $P_0(\theta)$. It will also be used to prove a pseudospectral estimate needed to compute the spectrum of $P_0(\theta)$:
\be
\sigma(P_0(\theta)) = \bigcup_{k\in\bN}\{k+ e^{-2i\varphi}r, r \ge 0\}
\ee
where $\varphi = \Im \theta$. Denote by $\vT= \bN$ the set of thresholds. Resonances of the KFP operator $P$ are defined as poles of the meromorphic continuation of the resolvent $R(z, \theta) =(P(\theta)  -z)^{-1} $, $\Im \theta >0$, across the positive semi-axis from the upper half-plane into the region in the lower half-plane defined by
\[
\{ z\in \bC, -2\varphi <\arg z \le 0\}\setminus \sigma(P_0(\theta)).
\]

As in quantum scattering, resonant states of the KFP operator are also metastable states. The real part of a resonance is related to the Kramers escape rate of the particle in Brownian motion and its inverse can be interpreted as the lifetime of the particle (see \cite{mel}). This is to be compared with resonances in quantum scattering where the lifetime is proportional to the inverse of the imaginary part.  \\

In certain problems such as resonances produced by tunneling effects, it is more convenient to use the method of analytic distortion outside a compact region (\cite{hu}). Here again the key point is  a global subelliptic estimate for the analytically distorted free KFP operator $\widetilde P_0(\theta)$. We prove that estimate (\ref{sub}) holds for $\widetilde P_0(\theta)$ and $n=1$ (Theorem \ref{th5.1}), but it is false for $n\ge 2$ (Proposition \ref{prop.nge2}). For this reason the definition of resonances based on analytic distortion is only valid  in dimension one. As expected, if both analytic dilation and analytic distortion are applicable, they give rise to the same set of resonances.
\\

The organization of this paper is as follows. In Section 2, we prove (\ref{sub}), which is the basic tool to define resonances of the KFP operator by analytic dilation, as well as $L^p$-subelliptic estimates, which are not used in the rest of this work.   In Section 3 we determine the spectrum of $P_0(\theta)$ and define resonances of the KFP operator $P$  as poles of meromorphic continuation  in $z$ of the resolvent $R(z, \theta)$ across the positive semi-axis into some region in the lower half-plane. It is proved that the poles in the positive semi-axis are eigenvalues and/or resonances of $P$ whose resonant states belong to weighted $L^2$-spaces and satisfy an outgoing radiation condition studied in \cite{w2}. The definition of resonances by analytic distortion is studied in Section 4. We only treat the one-dimensional case, due to the failure of a global subelliptic estimate for the distorted operator in higher dimensions.  \\

Notation. For $r\ge 0,s \in \bR$, we introduce  the weighted Sobolev space $\vH^{r,s}$ according to the subellipticity of $P_0$:
\[
\vH^{r,s} = \{u \in \vS'(\bR^{2n}), (1 + \w{D_v}^2 + |v|^2 + \w{D_x}^{\f 2 3})^{\f r 2 } \w{x}^s u \in L^2\}.
\]
For $ r <0, s\in \bR$, $\vH^{r, s}$ is the dual space of $\vH^{-r, -s}$. To simplify notation we denote  $\vH^s =\vH^{0,s}= L^2(\bR^{2n}, \w{x}^{2s} dxdv)$, $\vH = \vH^0= L^2(\bR^{2n})$.  $\vL(r,s, r',s')$ is the space of bounded linear operators from $\vH^{r,s}$ to $\vH^{r',s'}$ and $\vL(s,s') = \vL(0,s, 0,s')$. With this notation one has, for example,  $(P_0+1)^{-m}\in \vL(r, s, r+2m, s)$ for any $m\in \bN^*$ and $r,s \in \bR$. In this work, the multiplicity of an eigenvalue or a resonance refers to its algebraic multiplicity. The scalar product $\w{\cdot,\cdot}$ is linear in the second argument.

\sect{Global subelliptic estimates for complex vector fields}

For a closed linear operator $A$, let $\sigma (A)$ (resp., $\sigma_d(A)$, $\sigma_p(A)$) denote the spectrum (resp., the discrete spectrum, the point spectrum) of $A$.
There are several definitions of the essential spectrum of non-self-adjoint operators. We follow \cite{rs4} and define the essential spectrum of $A$ by
\[
\sigma_{\rm ess}(A) =\sigma(A)\setminus \sigma_d(A).
\]
With this definition, a relatively compact perturbation does not necessarily preserve the essential spectrum (\cite{rs4}).

 To determine  $\sigma(P_0(\theta))$  for $\theta \in \bC$, we regard $P_0(\theta)$ as a direct integral of the family of operators  $\{\widehat{P}_0(\xi, \theta) , \xi \in \bR^n\}$, where
\[
\widehat{P}_0(\xi, \theta) = \vF_{x\to \xi} P_0(\theta) \vF_{\xi \to x}^{-1} = H_0 + i v\cdot\xi_\theta.
\]
Here $\xi_\theta = e^{-\theta} \xi$. For $\zeta \in \bC^n$ with $\zeta =(\zeta_1, \cdots, \zeta_n)$, denote
 \[
\zeta^2 = \zeta \cdot \zeta =\sum_{j=1}^n \zeta_j^2.
\]
Then, using the properties of Hermite functions, one can check as in \cite{w1} that
\[
\sigma(\widehat{P}_0(\xi, \theta) )= \left\{k + \xi_\theta^2, k \in \bN\right \}.
\]
Each eigenvalue $k + \xi_\theta^2$ of $\widehat{P}_0(\xi, \theta)$ is semi-simple and the associated Riesz projection is given by
\be
\Pi_k^{\xi_\theta}f = \sum_{\alpha \in \bN^n, |\alpha |=k}\w{\overline{\psi_\alpha^{\xi_\theta}}, f} \psi_\alpha^{\xi_\theta}, \qquad f \in L^2(\bR^n_v).
\ee
Here $\{\psi_\alpha, \alpha \in \bN^n\}$ is an orthonormal basis of Hermite functions such that $ (-\Delta_v+\f{1}{4} |v|^2-\f{n}{2})\psi_\alpha = |\alpha| \psi_\alpha$ and
\[
\psi_\alpha^{\xi_\theta}(v) = \psi_\alpha(v+ 2i \xi_\theta).
\]

In order to define resonances of the KFP operator with a dilation-analytic potential, we need a global subelliptic estimate for $P_0(\theta) $ for $\theta \in \bC$. For $\Im \theta \neq 0$, $P_0(\theta) $ no longer satisfies H\"ormander's condition, because the vector field $e^{-\theta} v\cdot \partial_x$ is no longer real. A perturbation argument from the existing results for real vector fields  yields a global subelliptic estimate only when $|\Im \theta|$ is sufficiently small, but this kind of estimate is not convenient  for applications. The following parameter-dependent subelliptic estimate is important in this work. Its proof is based on an explicit formula for the kernel of the semigroup generated by $-e^{i\Im\theta}\widehat{P}_0(\xi, \theta)$ (Lemma \ref{lem2.2}), obtained by analytic continuation of the fundamental solution of the free KFP equation computed in \cite{wz}. \\

\begin{theorem} \label{th2.1}  For any $\phi_0\in ]0, \f \pi 2[$, there exists $C>0$ such that
\be \label{e2.2}
\| (s + H_0)f\| +  \w{\xi}^{\f 2 3 }\|f\|  + \| v\cdot \xi f\| \le C \|(\widehat{P}_0(\xi, \theta) +s) f\|,
\ee
 for $\theta \in \bC$ with $|\theta|\le \phi_0$,  $s \ge s_0 := \f n2$, $\xi\in \bR^n$,  and $f \in \vS(\bR_v^{n})$. Here $\|\cdot\|$ denotes the norm in $ L^2(\bR^n_v)$.
\end{theorem}
\medskip

Let $\theta = \mu + i \phi$ with $\mu, \phi\in \bR$, $|\phi| < \f \pi 2$, and let $\xi \neq 0$. Set $\omega = \xi/|\xi| \in \bS^{n-1}$ and $\rho = e^{-\mu}|\xi|$. Then
\be \label{e2.3}
\widehat{P}_0(\xi, \theta) = H_0 + i \zeta\cdot v = e^{-i\phi} G, \qquad \zeta := \rho e^{-i\phi}\omega, \quad G := e^{i\phi} H_0 + i \rho \, \omega\cdot v.
\ee
For $f\in \vS(\bR^n_v)$ one has $\Re \w{Gf, f} = \cos\phi \w{H_0 f, f} \ge 0$. The semigroup generated by $-G$ is $e^{-\tau \widehat P_0(\xi, \theta)}$ at the complex time $\tau = te^{i\phi}$, and its kernel can be obtained by analytic continuation of the fundamental solution of the free KFP operator computed in \cite{wz}.
For $\tau\in \bC$ with $\Re \tau >0$, let
\be \label{Ktau}
\vK_\tau(v,v') = \f{1}{(4\pi\sinh\tau)^{\f n2}} \exp\left(\f{n\tau}2 - \f{\coth \tau}{4}(v^2 + v'^2) + \f{{\rm cosech}\,\tau}{2}\, v\cdot v'\right), \qquad v, v'\in \bC^n,
\ee
where $v^2 = v\cdot v$ and the branch of $(\sinh\tau)^{\f n2}$ is the one which is positive for $\tau>0$. For $\tau >0$, $\vK_\tau$ is Mehler's kernel of $e^{-\tau H_0}$ (\cite[Lemma 2.1]{wz}). Set
\[
\sigma(\tau) = \tau - 2\tanh \f\tau 2, \qquad \Re \tau >0,
\]
and, for $\zeta\in\bC^n$,
\be \label{Ftau}
F_{\tau,\zeta}(v,v') = e^{-\tau\zeta^2}\vK_\tau(v+2i\zeta, v'+2i\zeta) = \vK_\tau(v,v') \exp\left(-i\tanh\f\tau 2\, (v+v')\cdot\zeta - \sigma(\tau)\zeta^2\right),
\ee
where the second equality is an elementary algebraic identity (for $\tau>0$ and $\zeta = \xi \in\bR^n$, this is the kernel of $e^{-\tau\widehat P_0(\xi, 0)}$ computed in \cite[Lemma 2.1]{wz}, with $\coth \tau -{\rm cosech}\,\tau = \tanh \f \tau 2$). Recall that $\psi_\alpha^\zeta(v) = \psi_\alpha(v+2i\zeta)$.

\begin{lemma}\label{lem2.2} Let $\rho \ge 0$, $\omega\in\bS^{n-1}$, $|\phi| < \f \pi 2$, $\zeta = \rho e^{-i\phi}\omega$ and $G$ be given by (\ref{e2.3}). For $t>0$, let $T(t)$ be the integral operator with kernel $F_{te^{i\phi},\zeta}(v,v')$, and $T(0) = 1$. Then:
\begin{enumerate}
\item[(a)] $T(t)$ is bounded on $L^2(\bR^n_v)$, $T(t+r) = T(t)T(r)$ and $\|T(t)\| \le 1$ for $t, r \ge 0$.
\item[(b)] For $t>0$ and $\tau = te^{i\phi}$,
\be \label{schur}
\|T(t)\| \le A(\tau)^n e^{-\rho^2 \widetilde E_\phi(t)},
\ee
where
\begin{align*}
A(\tau) & = \f{e^{\f t 2\cos\phi}}{\left(|\sinh\tau|\,\Re\coth\tau\right)^{\f 12}}, \\
\widetilde E_\phi(t) & = \Re\left(e^{-2i\phi}\sigma(\tau)\right) - \f{2\left(\Im \left(e^{-i\phi}\tanh\ds\f\tau2\right)\right)^2}{\Re\tanh\ds\f\tau 2}.
\end{align*}
\item[(c)] For $f\in \vS(\bR^n_v)$ and $t \ge 0$,
\[
T(t) f - f = - \int_0^t T(r) G f \, dr.
\]
\end{enumerate}
Moreover, for any $\phi_0\in ]0, \f \pi 2[$, there exist $t_1 \in ]0, 1]$ and $C_0\ge 1$ such that
\be \label{Esmall}
A(te^{i\phi})^n \le C_0, \qquad \widetilde E_\phi(t) \ge \f{\cos\phi}{24}\, t^3, \qquad 0< t\le t_1, \ |\phi|\le \phi_0.
\ee
\end{lemma}
\pf Let $t>0$ and $\tau = te^{i\phi}$, so that $\Re \tau = t\cos\phi>0$.
  Since $|e^{\tau}|>1$ and $w\mapsto \f{w+1}{w-1}$ maps $\{|w|>1\}$ onto $\{\Re w>0\}$, the numbers
\[
\Re\coth\tau, \qquad \Re\tanh\f\tau 2 = \Re(\coth\tau - {\rm cosech}\,\tau), \qquad \Re\coth\f\tau 2 = \Re(\coth\tau + {\rm cosech}\,\tau)
\]
are positive.  Put $a_\tau = \f14\Re\coth\tau$, $b_\tau = \f12\Re\,{\rm cosech}\,\tau$ and
\[
c = \Im\left(\tanh\f\tau2\, \zeta\right) = \rho\, \Im\left(e^{-i\phi}\tanh\f\tau2\right)\omega\in\bR^n.
\]
Then $2a_\tau >|b_\tau|$ and, for $v, v'\in\bR^n$,
\[
|F_{\tau,\zeta}(v,v')| = \f{e^{\f{nt}2\cos\phi - \rho^2\Re(e^{-2i\phi}\sigma(\tau))}}{(4\pi|\sinh\tau|)^{\f n2}}\, \exp\left(-a_\tau(|v|^2+|v'|^2) + b_\tau v\cdot v' + c\cdot(v+v')\right).
\]
A Gaussian integration gives
\[
\int_{\bR^n} \exp\left(-a_\tau|v'|^2 + (b_\tau v + c)\cdot v'\right) dv' = \left(\f\pi {a_\tau}\right)^{\f n2}\exp\left(\f{|b_\tau v+c|^2}{4a_\tau}\right),
\]
and an elementary computation shows that
\[
\sup_{v\in\bR^n}\left(-a_\tau|v|^2 + c\cdot v + \f{|b_\tau v+c|^2}{4a_\tau}\right) = \f{|c|^2}{2a_\tau-b_\tau} = \f{2|c|^2}{\Re\tanh\ds\f\tau 2},
\]
the quadratic form in $v$ having the negative coefficient $-\f{(2a_\tau-b_\tau)(2a_\tau+b_\tau)}{4a_\tau}$. Since $F_{\tau,\zeta}$ is symmetric in $(v,v')$, Schur's test yields
\[
\|T(t)\| \le \sup_{v}\int_{\bR^n}|F_{\tau,\zeta}(v,v')|\,dv' = A(\tau)^n e^{-\rho^2\widetilde E_\phi(t)},
\]
which is (\ref{schur}). In particular $T(t)$ is bounded on $L^2(\bR^n_v)$.

Next, note that
\be \label{eigen}
G\psi_\alpha^\zeta = e^{i\phi}\lambda_\alpha \psi_\alpha^\zeta, \qquad \lambda_\alpha = |\alpha| + \zeta^2.
\ee
For $\tau>0$ and $w\in\bR^n$, Mehler's formula gives $\int_{\bR^n}\vK_\tau(w,v')\psi_\alpha(v')\,dv' = e^{-\tau|\alpha|}\psi_\alpha(w)$. Since $\Re\coth\tau>0$, both sides are holomorphic in $(\tau, w)\in\{\Re\tau>0\}\times\bC^n$, so this identity holds for $\Re\tau>0$ and $w\in \bC^n$. The integrand $v'\mapsto \vK_\tau(v+2i\zeta, v')\psi_\alpha(v')$ is entire with Gaussian decay, uniformly in the strip between $\bR^n$ and $\bR^n + 2i\zeta$; shifting the contour $v'\to v'+2i\zeta$, we obtain
\be \label{Teigen}
T(t)\psi^\zeta_\alpha = e^{-\tau\zeta^2}\int_{\bR^n}\vK_\tau(v+2i\zeta, v')\psi_\alpha(v')\,dv' = e^{-\tau\lambda_\alpha}\psi_\alpha^\zeta, \qquad \tau = te^{i\phi}.
\ee

Let $\vP_\zeta$ be the linear span of $\{\psi_\alpha^\zeta, \alpha\in\bN^n\}$. $\vP_\zeta$ is equal to  the linear span of $\{e^{-i\Re\zeta\cdot v}\psi_\alpha(v-2\Im\zeta), \alpha\in\bN^n\}$, hence  is dense in $\vS$ and in $L^2(\bR^n_v)$.
By (\ref{Teigen}), $T(t)T(r) = T(t+r)$ on $\vP_\zeta$, hence on $L^2$ by density and boundedness. Let $f\in\vP_\zeta$. By (\ref{eigen}) and (\ref{Teigen}), $T(t)f\in\vP_\zeta\subset\vS$ and $\f{d}{dt}T(t)f = -GT(t)f$. Consequently
\[
\f{d}{dt}\|T(t)f\|^2 = -2\Re\w{GT(t)f, T(t)f} = -2\cos\phi\w{H_0T(t)f, T(t)f} \le 0,
\]
so $\|T(t)f\|\le\|f\|$, which proves (a). Next, (\ref{eigen}) and (\ref{Teigen}) give, for $f = \psi_\alpha^\zeta$,
\[
-\int_0^t T(r)Gf\,dr = -\int_0^t e^{i\phi}\lambda_\alpha e^{-re^{i\phi}\lambda_\alpha}\, dr\ \psi_\alpha^\zeta = T(t)f - f,
\]
hence (c) holds on $\vP_\zeta$. Since $\vP_\zeta$ is dense in $\vS$, $G$ is continuous from $\vS$ to $L^2$ and $\|T(r)\|\le 1$, (c) extends to $f\in\vS(\bR^n_v)$. \\

It remains to prove (\ref{Esmall}). One has
\[
|\sinh\tau|\,\Re\coth\tau = \f{\Re\left(\cosh\tau\,\overline{\sinh\tau}\right)}{|\sinh \tau|} \to \cos\phi \quad \mbox{ as } t\to 0^+,
\]
uniformly for $|\phi|\le\phi_0$, so $A(te^{i\phi})\to (\cos\phi)^{-\f12}$ and the first inequality in (\ref{Esmall}) holds for $t\le t_1$ with some $C_0$. Since
\[
\sigma(\tau) = \f{\tau^3}{12} + O(\tau^5), \qquad \tanh\f\tau 2 = \f\tau 2 - \f{\tau^3}{24} + O(\tau^5),
\]
one obtains
\[
\widetilde E_\phi(t) = \f{\cos\phi}{12}t^3 + O(t^5),
\]
uniformly for $|\phi|\le\phi_0$, and the second inequality in (\ref{Esmall}) holds for $0<t\le t_1$ if $t_1$ is small enough.
\ef

{\bf \noindent  Proof of Theorem \ref{th2.1}.} We prove the slightly more general estimate: for any $\phi_0\in ]0, \f \pi 2[$ there exists $C>0$ such that, with $G$ given by (\ref{e2.3}),
\be \label{e2.6}
\|(H_0+s)f\| + \w{\rho}^{\f 2 3}\|f\| + \rho\|\omega\cdot v f\| \le C \|(H_0 + i\rho e^{-i\phi}\omega\cdot v + s)f\|
\ee
for all $\rho\ge 0$, $\omega\in\bS^{n-1}$, $|\phi|\le\phi_0$, $s\ge s_0 = \f n2$ and $f\in \vS(\bR^n_v)$. Theorem \ref{th2.1} follows from (\ref{e2.6}) and (\ref{e2.3}). \\

Let $f\in\vS(\bR^n_v)$ and $u = (H_0 + i\rho e^{-i\phi}\omega\cdot v + s)f$, so that $(G + \lambda) f = e^{i\phi} u$ with $\lambda = s e^{i\phi}$.
By Lemma \ref{lem2.2} (c), $t\mapsto e^{-\lambda t} T(t) f$ is of class $C^1$ with derivative $-e^{-\lambda t}T(t)(G+\lambda)f$, and $\|e^{-\lambda t} T(t) f\| \le e^{-st\cos\phi}\|f\| \to 0$ as $t\to +\infty$. Hence
\[
f = \int_0^\infty e^{-\lambda t} T(t) (G+\lambda) f \, dt,
\]
and
\[
\|f\| \le J \|u\| \quad \mbox{ with } \quad J = \int_0^\infty e^{-st\cos\phi} \|T(t)\|\, dt.
\]
Clearly $J \le (s_0 \cos\phi_0)^{-1}$. Let $t_1$ be given by (\ref{Esmall}) and $c_1 = \f{\cos\phi_0}{24}$. By Lemma \ref{lem2.2}, $\|T(t)\| \le C_0e^{-c_1\rho^2 t^3}$ for $0< t\le t_1$, and $\|T(t)\| \le \|T(t-t_1)\|\,\|T(t_1)\| \le C_0e^{-c_1 \rho^2 t_1^3}$ for $t\ge t_1$. Therefore
\begin{align*}
J & \le C_0\int_0^\infty e^{-c_1\rho^2 t^3} dt + \f{C_0e^{-c_1\rho^2 t_1^3}}{s_0\cos\phi_0} \\
& \le C_0\Gamma\left(\f 4 3\right) (c_1 \rho^2)^{-\f 1 3} + \f{C_0e^{-c_1\rho^2 t_1^3}}{s_0\cos\phi_0} \le C \rho^{-\f 2 3}.
\end{align*}
Since also $J \le (s_0 \cos\phi_0)^{-1}$, we have proved that
\be \label{e2.7}
\w{\rho}^{\f 2 3} \|f\| \le C_1 \|u\|
\ee
uniformly in $\rho\ge 0$, $\omega$, $|\phi|\le\phi_0$ and $s\ge s_0$. \\

We now estimate the other terms in (\ref{e2.6}). Since $u = (H_0+s)f + \rho(\sin\phi + i\cos\phi)\, \omega\cdot v f$ and $[H_0, \omega\cdot v] = -2\omega\cdot\nabla_v$, one has
\begin{align*}
\|u\|^2 & = \|(H_0+s)f\|^2 + \rho^2\|\omega\cdot v f\|^2 \\
& \quad + 2\rho\sin\phi\, \Re \w{\omega\cdot vf, (H_0+s)f} + 2\rho\cos\phi \w{f, \omega\cdot D_v f},
\end{align*}
where $D_v = -i\nabla_v$. By the Cauchy-Schwarz inequality,
\[
2\rho|\sin\phi| \, |\w{\omega\cdot vf, (H_0+s)f}| \le \sin\phi_0 \left( \|(H_0+s)f\|^2 + \rho^2 \|\omega\cdot vf\|^2\right).
\]
For $s\ge \f n 2$, $\|\nabla_v f\|^2 \le \w{(H_0 + \f n 2) f, f} \le \|f\|\, \|(H_0+s)f\|$. Using (\ref{e2.7}) and $\rho \le \w{\rho}$, we get
\begin{align*}
2\rho |\w{f, \omega\cdot D_v f}| & \le 2\rho \|f\|^{\f 3 2} \|(H_0+s)f\|^{\f 12} \le 2C_1^{\f 3 2} \|u\|^{\f 3 2} \|(H_0+s)f\|^{\f 12} \\
& \le \f{1-\sin\phi_0}{2} \|(H_0+s)f\|^2 + C_2\|u\|^2
\end{align*}
by Young's inequality. Consequently,
\[
\f{1-\sin\phi_0}{2} \left(\|(H_0+s)f\|^2 + \rho^2\|\omega\cdot vf\|^2\right) \le (1 + C_2)\|u\|^2.
\]
Together with (\ref{e2.7}), this proves (\ref{e2.6}), hence (\ref{e2.2}).
\ef

\begin{rmk} \label{rmk2.4} The proof of (\ref{e2.6}) only uses the form of the operator $H_0 + i \zeta\cdot v$ with $\zeta = \rho e^{-i\phi}\omega$, $\rho\ge 0$, $\omega\in \bS^{n-1}$ and $|\phi|\le \phi_0$. In particular, (\ref{e2.2}) holds uniformly for $|\Re\theta|\le M$ and $|\Im \theta|\le \phi_0$, for any $M>0$ and $\phi_0\in ]0, \f \pi 2[$, and (\ref{e2.6}) will also be used in Section 4 for complex vector fields of the form $i v \eta$ with $\eta\in\bC$, $|\arg \eta| \le \phi_0$ ($n=1$).
\end{rmk}

Since $P_0(\theta)$ and $P_0(\theta')$ are unitarily equivalent when $\Im \theta = \Im \theta'$ (see also Remark \ref{rmk2.4}), Theorem \ref{th2.1} is true for any $\theta \in \bC$ with $|\Im \theta|< \f \pi 2$ and the estimate (\ref{e2.2}) is uniform for $|\Re \theta|\le M$ and $|\Im \theta |\le \phi_0$ (for fixed $M >0$ and $\phi_0 \in ]0, \f \pi 2[$). As a consequence of Theorem \ref{th2.1} (and the same estimate for $\widehat P_0(\xi, \theta)^*$), one has the following

\begin{cor} \label{cor2.2} Let $\vD= \{f\in L^2, P_0(\theta)f \in L^2\}$. For $\theta \in \bC$ with $|\Im \theta | < \f \pi 2$,  one has: $\vD = D(H_0) \cap D(v\cdot D_x)$.
 $P_0(\theta)$ defined on its maximal domain $\vD$ coincides with the closure of $P_0(\theta)$ on $\vS(\bR^{2n})$. Let $R_0(z,\theta) =(P_0(\theta)-z)^{-1}$ for $z\not\in\sigma(P_0(\theta))$. Then
\be\label{e2.12}
\|R_0(-s, \theta)\| \le \f C {\w{s}},
\ee
uniformly in $\theta \in \bC$ with $ | \theta | \le \phi_0 < \f \pi 2$ and $s \ge s_0(\phi_0)$.
\end{cor}

Under condition (\ref{ass1}), the global subelliptic estimate (\ref{e2.2}) shows that $W_\theta$ is a relatively compact perturbation of $P_0(\theta)$. Therefore $P(\theta)$ is closed with constant domain
$D(P(\theta)) =D(P_0(\theta))= \vD$. The family $\{P(\theta), \theta \in \bC, \, |\Im \theta | < a\}$ is a holomorphic family of type (A) in the sense of  Kato. Consequently the resolvent $R(z, \theta) =( P(\theta)-z)^{-1}$  is holomorphic in $(z,\theta)$ (\cite{k}).\\

The explicit kernel of Lemma \ref{lem2.2} also gives subelliptic estimates in $L^p$. For $1\le p\le \infty$, denote by $\|\cdot\|_p$ the norm in $L^p(\bR^{2n}_{x,v})$ and let $\w{D_x}^\nu = \vF^{-1}_{\xi\to x} \w{\xi}^\nu \vF_{x\to \xi}$; $\w{D_v}^\beta$ is defined in the same way.

\begin{theorem}\label{th-Lp} Let $1\le p<\infty$, $\phi_0\in ]0, \f\pi 2[$, $0\le \nu<\f 23$ and $0\le \beta< 2$. There exist $C, s_1>0$ such that
\be \label{eLp}
s\|f\|_p + s^{1-\f{3\nu}2}\|\w{D_x}^\nu f\|_p + s^{1-\f\beta 2}\left(\|\w{D_v}^\beta f\|_p + \|\w{v}^\beta f\|_p\right) \le C \|(P_0(\theta) + s)f\|_p
\ee
for $\theta\in\bC$ with $|\theta|\le\phi_0$, $s\ge s_1$ and $f\in\vS(\bR^{2n})$.
\end{theorem}
\pf Let $\theta = \mu + i\phi$, $t>0$ and $\tau = te^{i\phi}$. By Lemma \ref{lem2.2} applied with $\zeta = e^{-\theta}\xi$, $\xi\in\bR^n$, the semigroup $T(t)$ generated by $-e^{i\phi}P_0(\theta)$ in $L^2(\bR^{2n})$ has, for $0<t\le t_1$ with $t_1>0$ small, the kernel
\[
K_t(x,x',v,v') = \f{1}{(2\pi)^n}\int_{\bR^n} e^{i(x-x')\cdot\xi}F_{\tau, e^{-\theta}\xi}(v,v')\,d\xi.
\]
Set
\[
q = e^{-2\theta}\sigma(\tau), \qquad w = e^{-\theta}\tanh\f\tau 2, \qquad b = \Im w.
\]
Since $\sigma(\tau) = \f{\tau^3}{12}(1+O(t^2))$ and $\tanh\f\tau 2 = \f\tau 2 + O(t^3)$, one has $\arg q = \phi + O(t^2)$, $|q|\ge c\, t^3$, $|w| = O(t)$ and $b = O(t^3)$, uniformly for $|\theta|\le\phi_0$. In particular $\Re q>0$ if $t_1$ is small enough, so that the above integral converges absolutely, and (\ref{Ftau}) gives
\[
K_t(x,x',v,v') = \vK_\tau(v,v')\, g_q\left(x-x' - w(v+v')\right), \qquad g_q(y) = (4\pi q)^{-\f n2}e^{-\f{y^2}{4q}},
\]
$g_q$ being the heat kernel at the complex time $q$ (for $\theta = 0$, see \cite[Lemma 2.1]{wz}). We claim that, for $0<t\le t_1$, $|\theta|\le \phi_0$ and $1\le p\le \infty$,
\be \label{eTp}
\|T(t)\|_{p\to p} \le C, \quad \|\na_x T(t)\|_{p\to p}\le Ct^{-\f 32}, \quad \|\na_v^k T(t)\|_{p\to p}\le Ct^{-\f k2}, \quad \||v|^\beta T(t)\|_{p\to p}\le Ct^{-\f \beta 2},
\ee
where $k=1,2$ and $\beta\ge 0$. By Schur's test, it suffices to bound $\sup_{x,v}\int|K|\,dx'dv'$ and $\sup_{x',v'}\int|K|\,dxdv$ for the kernels $K$ of these operators. Set $z = v+v'$. Since $\Im (wz) = bz$, a Gaussian integration gives
\[
\int_{\bR^n} \left|g_q(y-wz)\right| dy = (\cos\arg q)^{-\f n2} e^{\kappa |z|^2}, \qquad \kappa = \f{b^2}{4|q|\cos\arg q} = O(t^3).
\]
With $a_\tau$, $b_\tau$ as in the proof of Lemma \ref{lem2.2}, one has
\[
|\vK_\tau(v,v')| = \f{e^{\f{nt}{2}\cos\phi}}{(4\pi|\sinh\tau|)^{\f n2}} \exp\left(-\f{2a_\tau-b_\tau}4 |v+v'|^2 - \f{2a_\tau + b_\tau}4 |v-v'|^2\right),
\]
where $2a_\tau - b_\tau = \f12\Re\tanh\f\tau 2 \ge c\,t$ and $2a_\tau + b_\tau = \f 12\Re\coth\f\tau 2 \ge \f c t$. Taking $t_1$ small, $\kappa \le \f{2a_\tau - b_\tau}8$, and the integration in $v-v'$ gives
\[
\int\!\!\int |K_t|\,dx'dv' \le C\,\left(|\sinh\tau|\,(2a_\tau+b_\tau)\right)^{-\f n2} \le C',
\]
since $|\sinh\tau|\,(2a_\tau+b_\tau) \to \cos\phi$ as $t\to 0^+$. By symmetry, this proves the first bound in (\ref{eTp}). The derivatives $\na_x$ and $\na_v$ of $K_t$ produce respectively the factors
\[
-\f{y-wz}{2q} \quad \mbox{ and } \quad \f{w(y-wz)}{2q} - \f14\tanh\f\tau 2\, (v+v') - \f14\coth\f\tau 2\, (v-v'), \qquad y = x-x',
\]
and $\na_v^2 K_t$ produces, in addition, the factor $-\f{w^2}{2q} - \f14\tanh\f\tau 2 - \f14\coth\f\tau 2 = O(t^{-1})$. Since $|q|^{-\f12} = O(t^{-\f32})$, $|v+v'|^\beta e^{-ct|v+v'|^2} = O(t^{-\f\beta 2})$, $|v-v'|^\beta e^{-c|v-v'|^2/t} = O(t^{\f\beta 2})$ and $|v|^\beta\le C_\beta(|v+v'|^\beta + |v-v'|^\beta)$, the other bounds in (\ref{eTp}) are obtained in the same way. \\

By (\ref{eTp}) and the semigroup property, $T(t)$ extends to $L^p$ and
\[
\|T(t)\|_{p\to p} \le Ce^{\omega_0 t}, \qquad \|AT(t)\|_{p\to p} \le C\, t^{-a}e^{\omega_0 t}, \qquad t>0,
\]
for some $\omega_0\ge 0$, where $(A, a) = (\na_x, \f 32)$, $(\na_v^k, \f k2)$ or $(|v|^\beta, \f\beta 2)$. By Lemma \ref{lem2.2} (c), $\|T(t)f -f\|_p\le Cte^{\omega_0 t}\|P_0(\theta)f\|_p$ for $f\in \vS$, so $\{T(t)\}$ is a $C_0$-semigroup on $L^p$, $p<\infty$. Using the formula
\[
\w{D}^\nu g = c_\nu \int_0^\infty \left(g - e^{-r}e^{r\Delta}g\right) r^{-1-\f\nu 2}\, dr, \qquad 0<\nu<2,
\]
in the variables $x$ or $v$, and $\|g - e^{r\Delta} g\|_p \le C\min\{\|g\|_p, r^{\f 12}\|\na g\|_p, r\|\na^2 g\|_p\}$, one obtains
\[
\|\w{D_x}^\nu g\|_p \le C\left(\|g\|_p + \|g\|_p^{1-\nu}\|\na_x g\|_p^{\nu}\right), \qquad \|\w{D_v}^\beta g\|_p \le C\left(\|g\|_p + \|g\|_p^{1-\f\beta 2}\|\na_v^2 g\|_p^{\f\beta 2}\right)
\]
for $0<\nu<1$ and $0<\beta<2$. Consequently $\|AT(t)\|_{p\to p} \le C (1+ t^{-a})e^{\omega_0 t}$ for $(A, a) = (\w{D_x}^\nu, \f{3\nu}2)$, $(\w{D_v}^\beta, \f\beta 2)$ and $(\w{v}^\beta, \f\beta 2)$. \\

Let $f\in\vS(\bR^{2n})$, $u = (P_0(\theta)+s)f$ and $\lambda = se^{i\phi}$. As in the proof of Theorem \ref{th2.1}, $f = e^{i\phi}\int_0^\infty e^{-\lambda t}T(t)u\, dt$, where the integral, convergent in $L^2$, also converges in $L^p$ if $s\cos\phi_0 \ge 2\omega_0$. For $A = 1$, $\w{D_x}^\nu$, $\w{D_v}^\beta$ or $\w{v}^\beta$ and the corresponding exponent $a\in [0, 1[$, $A$ is closed, and for $s\ge s_1 = \max\{1, \f{2\omega_0}{\cos\phi_0}\}$,
\[
\|Af\|_p \le C\int_0^\infty e^{-\f{st}{2}\cos\phi_0} \left(1+ t^{-a}\right)dt\ \|u\|_p \le C' s^{a-1}\|u\|_p.
\]
This proves (\ref{eLp}).
\ef

For $p=2$, the endpoint cases $\nu = \f 2 3$ and $\beta = 2$ of (\ref{eLp}) follow from Theorem \ref{th2.1}. For $p\neq 2$, they are not reached by the above semigroup argument, because $\|\w{D_x}^{\f 23}T(t)\|_{p\to p}$, $\|\na_v^2T(t)\|_{p\to p}$ and $\|\w{v}^2T(t)\|_{p\to p}$ are of order $t^{-1}$ as $t\to 0^+$.

\sect{Resonances by analytic dilation}

 We first define resonances of the KFP operator by analytic dilation. To begin with, we determine the spectrum of $P_0(\theta)$ by regarding it as a direct integral of the family of non-self-adjoint operators $\{\widehat{P}_0(\xi, \theta), \xi \in \bR^n\}$. For this purpose, we use the following

\begin{theorem}[\cite{w2}] \label{th3.2} Let $H$ be a separable Hilbert space, $X $ a nonempty open set of $\bR^n$ and $\vH = \{f : X \to H, \|f\| = (\int_X \|f(x)\|_H^2 dx)^{\f 1 2} < + \infty\}$ where $dx$ is the Lebesgue measure of $\bR^n$. Assume that  both $\{Q(x), x\in X\}$  and the adjoints  $\{Q(x)^*, x\in X\}$  are  strongly continuous families of closed, densely defined operators on $H$ with constant domains $D$ and $D^*$, respectively.
Assume in addition that  for each $z \in \bC \setminus \overline{ \cup_{x\in X} \sigma (Q(x))}$, one has
\be \label{e2.17}
\sup_{x\in X} \| (Q(x) -z)^{-1} \| < +\infty.
\ee
Let $Q$ be the closed, densely defined operator on $\vH$ such that for any $f$ in the domain of $Q$, one has
$f(x) \in D$ almost everywhere in $x$ and
\be
Qf = Q(x)f(x) \quad \mbox{ in } \vH.
\ee
Then one has
\be \label{eq2.18}
\sigma(Q) = \overline{ \cup_{x\in X} \sigma (Q(x))}.
\ee
\end{theorem}

Condition (\ref{e2.17}) is a kind of pseudospectral estimate for non-self-adjoint operators. It is trivial when $Q(x)$ is self-adjoint, because the norm of the resolvent of self-adjoint operators is bounded by  the inverse of the distance from $z$ to the spectrum. Without the condition (\ref{e2.17}), the equality (\ref{eq2.18}) may fail.

\begin{prop}\label{prop3.1} For $\theta \in \bC$ with $\varphi =\Im \theta \in ] -\f \pi 2, \f \pi 2[$, one has
\be
\sigma(P_0(\theta)) = \bigcup_{k\in\bN}\{k+ e^{-2i\varphi}r, r \ge 0\}.
\ee
\end{prop}
\pf
According to Theorem \ref{th3.2}, in order to prove that
\be \label{spectreP0}
\sigma(P_0(\theta) ) = \overline{\bigcup_{\xi\in \bR^n} \sigma(\widehat{P}_0(\xi, \theta))} = \bigcup_{k\in\bN}\{k+ e^{-2i \Im \theta}r, r \ge 0\},
\ee
the main task is  to check the  pseudospectral estimate:
\be \label{e2.19}
\sup_{\xi\in\bR^n} \|(\widehat{P}_0(\xi,\theta)-z)^{-1}\| \le C(z, \theta),
\ee
for $ z\not\in \overline{\cup_{\xi\in \bR^n} \sigma(\widehat{P}_0(\xi,\theta))}$. The other conditions of Theorem \ref{th3.2} are clearly satisfied. Estimate (\ref{e2.19}) is equivalent to the following lower bound
\be \label{pseudo}
\|(\widehat{P}_0(\xi,\theta)-z)f\| \ge c \|f\|, \qquad c>0,
\ee
uniformly in $\xi\in \bR^n$. We only need to prove (\ref{pseudo}) for $|\xi|$ large enough, because it holds true locally in $\xi$ for $ z\not\in \overline{\cup_{\xi\in \bR^n} \sigma(\widehat{P}_0(\xi,\theta))}$. Applying (\ref{e2.2}) with a fixed $s \ge s_0$, one has
\begin{align*}
\lefteqn{ \|(\widehat{P}_0(\xi,\theta)-z)f\|} \\
& \ge \|(\widehat{P}_0(\xi,\theta)+ s)f\| -(s + |z|) \|f\| \\
&\ge \left( C^{-1}  \w{\xi}^{\f 2 3} -(s +|z|)\right) \|f\| \ge (2C)^{-1} \w{\xi}^{\f 2 3} \|f\|
\end{align*}
for $|\xi|\ge \xi_0$ with $\xi_0>0$ appropriately large. This proves the estimate (\ref{pseudo}). Proposition \ref{prop3.1} follows from Theorem \ref{th3.2}. \ef

\begin{prop}\label{prop2.3} Assume the condition (\ref{ass1}).
 Then for $\theta \in \bC$ with $|\Im \theta | < a$,  $P(\theta)$ defined on $\vD = D(H_0) \cap D(v\cdot D_x)$ is closed and
 \be \label{essP}
 \sigma_{\rm ess}(P(\theta)) = \bigcup_{k\in\bN}\{k+ e^{-2i\varphi}r, r \ge 0\}, \qquad \varphi = \Im \theta.
 \ee
\end{prop}
\pf By Theorem \ref{th2.1} and condition (\ref{ass1}), $W_\theta$ is relatively compact with respect to $P_0(\theta)$.
But with our definition of essential spectrum, this property alone is not sufficient to conclude (\ref{essP}). Proposition \ref{prop3.1} shows that $\sigma_{\rm ess}(P_0(\theta) ) = \sigma(P_0(\theta))$ has empty interior and that $\bC \setminus \sigma (P_0(\theta))$ is connected and contains $]-\infty, -s_0]$. Theorem \ref{th2.1} also implies that
\[
\|W_\theta R_0(-s, \theta) \| \le C \w{s}^{-\f 1 2}, \qquad s \ge s_0.
\]
Whence $(P(\theta)+s)^{-1}$ exists for $s \ge s_1$ if $s_1$ is appropriately large. With these properties, we can now apply Lemma 3 in Section XIII.4 of \cite{rs4} on perturbation of essential spectrum  to conclude that
\[
 \sigma_{\rm ess}(P(\theta)) =  \sigma_{\rm ess}(P_0(\theta)).
\]
This proves  (\ref{essP}).
\ef

In the following, we use the convention that for $z\in \bC^*$, $\arg z $ varies continuously in $]-\f \pi 2, \f{3\pi}2[$. For $\theta \in \bC$ with $\Im \theta >0$, denote
\be
 \Omega_-(\theta) = \{ z\in \bC^*, -2 \Im \theta < \arg z \le 0\} \setminus \left( \bigcup_{k\in\bN}\{k+ e^{-2i\Im \theta }r, r \ge 0\}\right).
\ee
   Let $\bC_0 =\{z\in \bC^*, 0<\arg z < \f {3\pi} 2 \}$. For $0\le \Im \theta <a$,  the spectrum of $P(\theta)$ in $\bC_0$ is discrete and independent of $\theta$; it consists of the eigenvalues of $P$ in $\bC_0$.
Denote $R_0(z, \theta) = ( P_0(\theta)-z)^{-1}$  and  $R(z, \theta) = ( P(\theta)-z)^{-1}$. By Proposition \ref{prop2.3}, $R(z, \theta)$ is an  operator-valued meromorphic function for $z\in \bC_0$.   $R_0(z, \theta)W_\theta$ is a holomorphic family of compact operators for $z\in \bC_0\cup \Omega_-(\theta)$ and  $\theta \in \bC$ with $0< \Im \theta < a$, and
\[
\|R_0(-s, \theta) W_\theta\| \le C\w{s}^{-\f 1 2}
\]
 for $s\ge s_0$. Hence there exists $s_0'>0$ such that $1+ R_0(-s, \theta) W_\theta$ is invertible for all $s\ge s'_0$ and
\be \label{resol-dilation}
R(z, \theta)  = \left( 1 + R_0(z, \theta) W_\theta \right )^{-1} R_0(z, \theta)
\ee
holds for $z=-s$ and $s\ge s'_0$.   The analytic Fredholm theorem (Theorem XIII.13, \cite{rs4}) shows that
  $ (1+ R_0(z, \theta) W_\theta)^{-1}$  is  meromorphic in $z\in \bC_0\cup \Omega_-(\theta)$. Equation (\ref{resol-dilation}) gives a meromorphic continuation  of $ R(z, \theta)$ from $\bC_0$ across the positive real semi-axis into $\Omega_-(\theta)$.
    \\

    The poles of $R(z, \theta)$ are independent of $\theta$  in the sense that if $\theta$ and $\theta'$ are two complex numbers such that $0 < \Im \theta < \Im \theta' < a$,  then $R(z, \theta')$ has the same poles as $R(z, \theta)$ in $\Omega_-(\theta)$ with the same algebraic multiplicities (see \cite{ac}).  As noted above, poles in $\bC_0$  of $ R(z, \theta)$  for $ 0 \le \Im \theta< a$ are complex eigenvalues of $P$.  \\

  In summary, we have proved the following

\begin{theorem}\label{th3.3}  Assume condition (\ref{ass1}) for some $a \in ]0, \f \pi 4]$. For $ 0< \Im \theta <a$, the operator-valued meromorphic function $R(z, \theta)$ defined for $z$ in $\bC_0$ can be meromorphically extended across the positive real semi-axis into $\Omega_-(\theta)$.
\end{theorem}

\begin{defn}\label{defn3.4}
 The set of resonances of $P$ is defined as
\be \label{resonances}
{\rm Res } (P) = \bigcup_{\theta \in \bC^*,  0< \Im \theta < a} \{ z \in \Omega_-(\theta), \mbox{ $z$ is a pole of
$R(z, \theta)$ }  \}.
\ee
\end{defn}

Note that the domain of definition of resonances for the KFP operator is analogous to that for the Schr\"odinger operator with a constant magnetic field where thresholds are Landau levels (\cite{bbr}).\\

$R(z, \theta)$ is meromorphic in the region $\bC_0 \cup \Omega_-(\theta)$. Poles of $R(z, \theta)$ in  $\bC_0$ are discrete complex eigenvalues of $P$; they have non-negative real parts, since $\Re \w{Pf, f} = \w{H_0 f, f} \ge 0$. The above definition of resonances includes poles of $R(z, \theta)$ in $\bR^*_+\setminus\vT $. For the Schr\"odinger operator with a complex short-range potential, positive eigenvalues are absent, but positive resonances may exist, and
positive poles of the meromorphic extension of the resolvent coincide with positive resonances whose resonant states belong to weighted $L^2$-spaces and satisfy the outgoing Sommerfeld radiation condition. See \cite[Theorem 5.3]{w3}.  For the KFP operator,   we prove the following

\begin{prop}\label{prop3.4} Assume that $V$ is long-range: $V(x) = O(\w{x}^{-\ep_0})$,  and satisfies (\ref{ass1}) with $\rho=1+\ep_0$ for some $\ep_0 >0$. Let $\lambda \in \bR_+ \setminus \vT$.  Then $\lambda \in  {\rm Res } (P)  $ if and only if
the equation $(1+ R_0(\lambda+i0) W)u =0$ has a nontrivial solution in $\vH^{-s}$ for  some (equivalently, for every) $s\in\, ]\f 12, \f{\rho}2[$. Here $R_0(z) = (P_0 -z)^{-1}$ and
\[
R_0(\lambda\pm i0) = \lim_{\pm \Im z >0, \, z\to \lambda} R_0(z)
\]
 in $\vL(s, -s)$, $s > \f 12$.
\end{prop}
\pf First note that real resonances are poles of $R(z, \theta)$ for any $\theta$ with $0<\Im \theta <a$. Let $\lambda \in \bR_+ \setminus \vT$. By the limiting absorption principle (\cite{w2}) for $P_0$, for $s> \f 12 $ there exists some constant $\delta>0$ such that
\be \label{LAP-R0}
R_0(z) - R_0(z') = O(|z-z'|^\delta) \qquad \mbox{ in } \vL(-1,s,1,-s),
\ee
for $z, z'$ near $\lambda$ with $\Im z>0$ and $ \Im z' >0$.
For $z$ near $\lambda$ with $\Im z >0$ and $|\theta|$ small with $\Im \theta \ge 0$, let $z' = e^{\theta} z$. Then
$\Im z'>0$ and, since $P_0(\theta) - z = e^{-\theta}(e^{\theta}H_0 + v\cdot\na_x - z')$,
\[
R_0(z, \theta) = e^{\theta} (e^{\theta }H_0 + v\cdot\na_x -z')^{-1}.
\]
Following the proof of (\ref{LAP-R0}) given in \cite{w2}, one can show that
\be\label{e3.41}
R_0(z, \theta)   = R_0(z') + O(|\theta|^\delta)
\ee
in $\vL(-1,s,1,-s)$, $s > \f 12$. Since $|z-z'| = O(|\theta|)$, it follows from (\ref{LAP-R0}) and (\ref{e3.41}) that
\begin{align*}
R_0(z) - R_0(z, \theta)
& = (R_0(z) - R_0(z') ) + ( R_0(z') -  R_0(z, \theta) )\\
& = O(|\theta|^\delta)
\end{align*}
in $\vL(-1,s,1,-s)$, $s > \f 12$, uniformly in $z$ near $\lambda$ with $\Im z>0$. For $V$ satisfying (\ref{ass1}) with $\rho >1$, one has
\be
R_0(z)W - R_0(z, \theta)W_\theta = O(|\theta|^\delta) \qquad \mbox{ in } \vL(\vH^{-s}), \quad  \f 12< s < \f\rho 2, \label{e3.41b}
\ee
uniformly in $z$ near $\lambda$ with $\Im z >0$ and $ \theta$ near $0$ with $\Im \theta \ge 0$. Fix $s\in ]\f 12, \f\rho 2[$ and set $K(0) = R_0(\lambda+i0)W$ and $K(\theta) = R_0(\lambda, \theta)W_\theta$ for $0<\Im\theta<a$. They are compact operators on $\vH^{-s}$ and $K(\theta)$ is holomorphic in $\theta$. Taking the limit $z \to \lambda+i0$ in (\ref{e3.41b}), one obtains
\be \label{e3.42}
\|K(0) - K(\theta)\|_{\vL(\vH^{-s})} \le C|\theta|^\delta
\ee
for $|\theta|$ small and $\Im\theta>0$.

By the dilation analyticity of $V$,  the eigenvalues of $K(\theta)$ are independent of $\theta$ with  $\Im \theta>0$ \cite{ac}.
 Let $\Gamma = \{\mu\in \bC, |\mu+1| = \ep\}$, where $\ep>0$ is such that $\sigma(K(0))\cap\{|\mu+1|\le\ep\}\subset\{-1\}$. By (\ref{e3.42}), there exists $r>0$ such that $\Gamma$ is contained in the resolvent set of $K(\theta)$ for $\theta\in D_+ = \{\theta\in\bC, |\theta|<r, \Im\theta>0\}$. Define
 \[
 \Pi_{-1} = \f{i}{2\pi} \int_{|\mu +1|=\ep} (K(0) -\mu)^{-1} d \mu, \qquad \Pi_{-1}(\theta) = \f{i}{2\pi} \int_{|\mu +1|=\ep} (K(\theta) -\mu)^{-1} d \mu.
 \]
 Estimate (\ref{e3.42}) shows that $\| \Pi_{-1}(\theta) -\Pi_{-1}\| = O(|\theta|^\delta)$ for $\theta \in D_+$. Therefore $K(\theta)$ has the same number of eigenvalues as $K(0)$, counted with multiplicity, inside $\Gamma$. In particular, if $-1$ is not an eigenvalue of $K(0)$, then $-1$ is not an eigenvalue of $K(\theta)$ for $\theta\in D_+$ (with $r>0$ appropriately small). Conversely, if $-1$ is an eigenvalue of $K(0)$ with multiplicity $m$, then by finite-dimensional perturbation theory of eigenvalues \cite{k}, $K(\theta)$ has $m$ eigenvalues, $\mu_j$, counted with multiplicity, satisfying
\[
\mu_j \to -1 \quad \mbox{ as } \theta\to 0, \ \theta \in D_+.
\]
 Since the $\mu_j$ are independent of $\theta$, one obtains $\mu_j= -1$, hence $-1$ is an eigenvalue of $K(\theta)$ with multiplicity $m$. Consequently,  $1+K(0)$ is injective in $\vH^{-s}$ if and only if $1+K(\theta_0)$ is injective in $\vH^{-s}$ for $\theta_0\in D_+$.

If $1+K(0)$ is not injective, there exists $u\in \vH^{-s}$, $u\neq 0$, such that $(1+K(\theta_0))u=0$. Since $\lambda \not \in \sigma(P_0(\theta_0))$, a bootstrap argument shows that $u \in \vH$  and $(P(\theta_0) -\lambda)u=0$, hence $R(z, \theta_0)$ has a pole at $\lambda$. Conversely, if $\lambda\in{\rm Res}(P)$, then $\lambda$ is an eigenvalue of $P(\theta_0)$ for some $\theta_0\in D_+$, hence $1+K(\theta_0)$ is not injective on $\vH$, which implies that
$1+K(0)$ is not injective on $\vH^{-s}$.
\ef

A nontrivial solution $u$ of $(1+ R_0(\lambda+i0) W)u =0$ with $u\in \vH^{-s}$ is either an eigenfunction or a resonant state satisfying the outgoing radiation condition discussed in \cite{w2}. Proposition \ref{prop3.4} identifies poles of $R(z, \theta)$ in $\bR_+\setminus \vT$ with embedded eigenvalues and/or outgoing resonances of $P$.

\sect{Resonances by analytic distortion in dimension one}

To study resonances induced by the tunneling effect in the Kramers problem, it is more convenient to define resonances by using the  method of analytic distortion outside a compact domain in $\bR^n_x$ (see \cite{hu} for Schr\"odinger operators). Due to the failure of subellipticity  in higher dimensions, we only give a definition in the one-dimensional case.
\\

 Let $r_0>1$  and $\chi  \in C^\infty(\bR)$ be real-valued such that $ 0 \le \chi \le 1$ and $\chi (r) =0$ if $ r \le r_0$,
 $\chi (r) =1$ if $ r \ge C r_0 $, $C>1$, and $\chi '(r) \ge 0$ such that
 \[
 \f{1}{\ln C} < \sup_{r} |r\chi'(r)| < \f{2}{\ln C} .
 \]
 One can check that functions satisfying these requirements exist. Consider the map $ (x, v) \to (F_\theta(x), v)$ in $\bR^{2n}$ where
 \be \label{e5.1}
 F_\theta(x) = e^{\theta \chi (|x|)}x, \quad x \in \bR^n.
 \ee
 By the choice  of $\chi$ with $C> e^{\f \pi 2}$,  $x\to F_\theta(x)$ is a global diffeomorphism in $\bR^n$ for $\theta \in \bR$ with $|\theta| \le \f \pi 4$.
Set
\be \label{e5.2}
\widetilde U_\theta f (x,v) = |DF_\theta(x)|^{\f 12} f(F_\theta(x),v), \quad f \in L^2(\bR_{x,v}^{2n}),
\ee
where $DF_\theta (x)$ is the Jacobian matrix and $ |DF_\theta(x)|$  the Jacobian determinant of  the change of variables: $x \to F_\theta(x)$.
One has
\[
DF_\theta (x) = e^{\theta \chi(|x|)} \left(1 +  \theta \chi'(|x|) \f{x \otimes x}{|x|}\right)
\]
and $ A^\theta(x)= (DF_\theta(x))^{-1}$ is given by
\[
A^\theta(x) = e^{-\theta \chi(|x|)} \left(1 -\f{\theta \chi'(|x|)}{(1+\theta \chi'(|x|)|x|)} \f{x \otimes x}{|x|}\right).
\]
For $\theta \in \bR$ with  $|\theta| \le \f \pi 4$,
 define the distorted operator $\widetilde {P}_0(\theta)$ by
\[
\widetilde{P}_0(\theta) = \widetilde U_\theta P_0 \widetilde U_\theta^{-1}.
\]
One has
\be \label{e5.3}
\widetilde P_0(\theta) = H_0 + v\cdot \nabla_x^\theta  -\f 1 {2 |DF_\theta(x)|} v\cdot A^\theta(x)(\nabla |DF_\theta(x)|)
\ee
where
\[
 \nabla^\theta_x  =  A^\theta(x) \nabla_x.
\]
By (\ref{e5.3}) $\widetilde P_0(\theta)$ is well-defined for $\theta \in \bC$ with $|\theta|\le \f \pi  4$.
Denote
\be
{\widetilde Q}_0(\eta^\theta) = H_0 + i v\cdot \eta^\theta,
\ee
where
\[
\eta^\theta = A^\theta(x) \xi.
\]

\begin{theorem} \label{th5.1}  Let $n=1$. Then there exist some positive constants $ C$ and $s_0$ such that
\be \label{e5.4}
\| (s + H_0)f\| +  \|v\nabla_x f\| + \|\w{D_x}^{\f 2 3 }f\|  \le C \|({\widetilde P}_0(\theta) +s) f\|,
\ee
  uniformly in $\theta \in \bC$ with $|\theta|\le \f{\pi}4$,  $s \ge s_0$, and $f \in \vS(\bR_{x,v}^{2})$.
\end{theorem}
\pf
 First,  we show that
\be \label{e5.5}
\| (H_0+s)f\| + \|v \eta^\theta f\| + \w{\xi}^{\f 2 3}\|f\| \le C\|({\widetilde Q}_0(\eta^\theta) +s)f\|
\ee
for $\xi \in \bR$, $s \ge s_0$, $\theta \in \bC$ with $ |\theta | \le \f \pi 4$ and  $f\in \vS(\bR_v)$, uniformly in $x\in\bR$. Here $\|\cdot\|$ is the norm in $L^2(\bR_v)$.
In fact, for $n=1$, $A^\theta(x)$ is a complex number.
Since $0\le \chi\le 1$, $0\le |x|\chi'(|x|) \le \f 2{\ln C}$ and $|\theta|\le\f\pi 4$, one has, if the constant $C$ in $\chi$ is taken large enough,
\[
c \le |A^\theta(x)| \le c^{-1}, \qquad |\arg A^\theta(x)| \le \f \pi 3,
\]
for some $c>0$ independent of $x$ and $\theta$. Write $\eta^\theta = A^\theta(x)\xi = \rho e^{-i\phi}\omega$ with $\rho = |A^\theta(x)|\,|\xi|$, $\omega = \sgn \xi$ and $\phi = -\arg A^\theta(x)$, $|\phi|\le \f \pi 3$. Then ${\widetilde Q}_0(\eta^\theta) = H_0 + i\rho e^{-i\phi}\omega v$, and (\ref{e5.5}) follows from the estimate (\ref{e2.6}) established in the proof of Theorem \ref{th2.1} (see Remark \ref{rmk2.4}), applied with $\phi_0 = \f\pi 3$, since $\w{\rho}\ge c\w{\xi}$ and $\rho\|\omega v f\| = \|v\eta^\theta f\|$.\\

 Next, to derive the global subelliptic estimate (\ref{e5.4}), we now construct a parametrix for $\widetilde{ P}_0(\theta)+s$ using pseudodifferential operators with operator-valued symbols (\cite{maso}).  By (\ref{e5.5}) (and its adjoint version), $(\widetilde{Q}_0(\eta^\theta) +s)^{-1}$ exists for $s \ge s_0$. Set
 \be
 a(x, \xi , s) = (\widetilde{Q}_0(\eta^\theta) +s)^{-1}, \quad (x, \xi) \in \bR^2.
 \ee
 By (\ref{e5.5}),
 \[
 \|(H_0 + s+ |v\xi| + \w{\xi}^{\f 23}) a(x, \xi , s)\|\le C
 \]
 uniformly in $x$ and $\xi$.
Since $n=1$, one has  for $\alpha \ge 1$
 \be
 |\partial^{\alpha}_x v \eta^\theta | \le C_\alpha |\theta| (1+ |v\xi|)
 \ee
 and
 \be
  |\partial^{\alpha}_\xi v  \eta^\theta | \le C_\alpha |v|.
  \ee
 This gives for $\alpha, \beta \in \bN$,
 \be\label{symbol}
 \| \partial^\alpha_x \partial^\beta_\xi a(x, \xi, s)\| \le C_{\alpha \beta} \|a(x, \xi , s)\| \w{\xi}^{-\f 1 3|\beta|}.
\ee
Note that $ a(x, \xi , s)  = a_{\infty}(\xi, s) $  for $|x| \ge C r_0$, where
\[
a_{\infty}(\xi, s) = (\widetilde{Q}_0(e^{-\theta} \xi) +s)^{-1}.
\]
Let $a(x, D_x,s)$  denote the pseudodifferential operator with operator-valued symbol $a$ defined by
\[
(a(x, D_x, s)f)(x, v) = \f{1}{2\pi} \int\int_{\bR^{2}} e^{i (x-y)\cdot\xi} \left(a(x, \xi, s)f(y, \cdot)\right)(v)  d\xi dy
\]
for $f\in \vS(\bR^{2}_{x,v})$. Then
\be \label{e5.30}
(H_0+s)a(x, D_x, s), \quad \w{D_x}^{\f 2 3} a(x, D_x, s) \quad \mbox{  and } v\cdot \nabla_x a(x, D_x, s)
\ee
are uniformly bounded for $s$ large (\cite{maso}).
A direct computation gives
\be
(\widetilde{ P}_0(\theta)+ s)a(x, D_x, s) = 1 + b(x, D_x,s) + O(\w{s}^{-\f 12}),
\ee
where
\[
b(x, \xi, s) =  \nabla_\xi(v\cdot \eta^\theta ) \nabla_x a(x, \xi, s),
\]
and the term $O(\w{s}^{-\f 12})$ comes from the zero-order term in (\ref{e5.3}), because $\|v a(x,\xi,s)\| \le C\w{s}^{-\f 12}$.
The operator-valued symbol $b(x, \xi, s) $ satisfies the estimates
\beas
\|b(x, \xi, s)\| & \le & C|\theta| \|   v  a(x, \xi, s)\| \, \|  v\xi  a(x, \xi, s) \|\\
& \le & C |\theta| \|a(x, \xi, s)\|^{\f 1 2} \le C'|\theta|\w{s}^{-\f 1 2}
\eeas
and
\[
 \|\partial^\alpha_x \partial^\beta_\xi b(x, \xi, s)\| \le C_{\alpha \beta} |\theta| \w{s}^{-\f 1 2} \w{\xi}^{ -\f 1 3|\beta|}, \quad \alpha, \beta \in \bN,
\]
uniformly in $x, \xi $ and $s$.
It follows that $b(x,D_x,s)$ is bounded and
\[
\|b(x, D_x,s)\| \le C |\theta| \w{s}^{-\f 1 2}.
\]
This shows  that there exists $s_0' \ge s_0$ such that $\widetilde{ P}_0(\theta)+ s$ has a right inverse for $s \ge s'_0$ and
\[
(\widetilde{ P}_0(\theta)+ s) a(x, D_x, s)( 1 + O(\w{s}^{-\f 1 2})) = 1.
\]
Similarly, one can construct a left inverse of  $\widetilde{ P}_0(\theta)+ s$  for $ s$ large. Therefore  $\widetilde{ P}_0(\theta)+ s$ is invertible  for $s \ge s'_0$, with $s'_0>1$ appropriately large, and
\[
(\widetilde{ P}_0(\theta)+ s)^{-1} =  a(x, D_x, s)( 1 + O(\w{s}^{-\f 1 2})).
\]
Estimate (\ref{e5.4}) for $s\ge s_0'$ follows from (\ref{e5.30}).
\ef

When $n \ge 2$, $|\partial^\alpha_x (v \cdot \eta^\theta)|$ cannot be bounded by $|v \cdot\xi|$ and $\w{\xi}^{\f 2 3}$.
This explains why our approach fails for dimensions $n \ge 2$.  One may wonder if one can use other methods,  like Kohn's method (\cite{hln}) or the method of multipliers (\cite{hl}),  to prove a global subelliptic estimate for $\widetilde P_0(\theta)$ when $n\ge 2$. The following result gives a negative answer to this question.

\begin{prop}\label{prop.nge2} Let $n \ge 2$ and $\f 3 4 < \sigma \le 1$. Then for $|\theta| < \f \pi 4$ and $\Im \theta \neq 0$, $H_0^\sigma$ is not relatively bounded with respect to $\widetilde{ P}_0(\theta)$, that is, there is no constant $C_0$ such that
\be \label{nrb}
\|H_0^\sigma u\| \le C_0 ( \| \widetilde{ P}_0(\theta) u\| + \| u\|), \qquad \forall u \in \vS(\bR^{2n}).
\ee
\end{prop}
\pf We argue by contradiction and construct a sequence of functions $u_j \in C_0^\infty(\bR^{2n})$ invalidating (\ref{nrb}) for any $C_0>0$.
To simplify notation, let $\theta= i\phi$ with $0<\phi <\f \pi 4$. When $n \ge 2$, $A^\theta(x)$ has two eigenvalues
\[
a(x) =\f{e^{-i\phi \chi(|x|)}(1-i \kappa)}{1 + \kappa^2}, \quad b(x) =e^{-i\phi \chi(|x|)}
\]
where $\kappa (x) = \phi  r\chi'(r)$, $r =|x|$. Denote
\[
a(x) = a_1(x) + i a_2(x), \quad b(x) = b_1(x) + i b_2(x)
\]
with $a_j(x), b_j(x)$ real-valued.  Fix a point $x_0 $ such that  $\chi'(|x_0|) >0$. Denote $A^\theta = A^\theta(x_0)$, and similarly for the functions introduced above evaluated at $x_0$. Then $\kappa, a_1, b_1$ are positive and
\[
a_2 = - \f{\sin(\phi \chi) + \kappa \cos(\phi \chi)}{1+\kappa^2}<0, \quad b_2 = -\sin(\phi \chi) <0.
\]
 Let $\omega \in \bS^{n-1}$ with $\omega \perp \hat x_0$.
For $\xi = s \hat x_0 + t \omega$ and $v = \lambda \hat x_0 + \mu \omega$ with $s, t, \lambda, \mu  \in \bR$, one has
\[
v\cdot A^\theta \xi = a s \lambda + b t \mu=  (a_1 s \lambda + b_1 t  \mu ) + i (a_2 s \lambda + b_2 t \mu).
\]
Consequently, if $v = \lambda (\hat x_0 + \omega)$ and $\xi= \tau (b_1 \hat x_0 -a_1 \omega)$ with $\tau \in \bR$, then
\beas
\Re (v\cdot A^\theta \xi )  &= &0, \\
\Im (v\cdot A^\theta \xi) & = &\lambda \tau ( a_2 b_1- a_1 b_2) = - \lambda \tau \f{\kappa}{1+\kappa^2}, \\
\f{|v|^2} 4 - \f n 2 + i v\cdot A^\theta \xi & = &\f{\lambda^2} 2 - \f n 2 + \f{\lambda \tau \kappa}{1+\kappa^2}.
\eeas
This shows that for $v_j = j (\hat x_0 + \omega)$ and $ \xi_j = \tau_j (b_1 \hat x_0 - a_1 \omega)$ with
\[
\tau_j = -\f{(1+\kappa^2) (j^2-n)}{2j \kappa},
\]
one has
\be
\f{|v_j|^2} 4 - \f n 2 + i v_j\cdot A^\theta \xi_j = 0, \qquad  j \in \bN^*. \label{charact}
\ee
Let $f, g \in C_0^\infty(\bR^n)$ with $\|f\| =\|g\|=1$. Define
\[
u_j(x, v) = j^{\f{5n}{12}} e^{i(x-x_0)\cdot \xi_j} f(j^{\f 1 2}(x-x_0)) g(j^{\f 1 3 } ( v-v_j)).
\]
Then $u_j \in C_0^\infty(\bR^{2n})$ and $\|u_j\|=1$. One has
\beas
\|\Delta_v u_j \|  &= & O(j^{\f 2 3}), \\
\||v|^{2\sigma } u_j\| &\ge & \| |v_j +  O(j ^{-\f 1 3})|^{2 \sigma} u_j\| \ge j^{2\sigma}
\eeas
for $j \ge j_0$ with $j_0 $ appropriately large. Note that $(H_0 + n)^{-\sigma}$ is a pseudodifferential operator with principal symbol
$(|\eta|^2 +\f 1 4 |v|^2 +\f n 2)^{-\sigma}$  (\cite{hor}).  Therefore $|v|^{2\sigma} (H_0 + n)^{-\sigma}$ is bounded: $\exists C >0$ such that
\[
\||v|^{2\sigma} u\| \le C \| (H_0 + n)^\sigma u\| \le C (\| H_0^\sigma u\| + n^\sigma \|u\|), \qquad u \in \vS(\bR^{2n}).
\]
It follows that there exists $c_0>0$ such that
\be \label{lowerb}
\|H_0^\sigma u_j\| \ge c_0 j^{2\sigma}
\ee
for $j \ge j_1$ with $j_1$ appropriately large. On the other hand,
\[
\| \widetilde{ P}_0(\theta) u_j\| \le O(j) + \left \|\left(\f{|v|^2} 4 - \f n 2 + v\cdot A^\theta(x) \nabla_x \right )u_j \right \|,
\]
where $O(j)$ comes from $\Delta_v u_j$ and from the zero-order term in (\ref{e5.3}). On the support of $u_j$,  $A^\theta(x) - A^\theta(x_0) = O(j^{-\f 1 2})$, $v-v_j = O(j^{-\f 13})$ and $|v| + |\xi_j| = O(j)$. Using (\ref{charact}), one can evaluate
\beas
\lefteqn{ \left(\f{|v|^2} 4 - \f n 2 + v\cdot A^\theta(x) \nabla_x\right )u_j } \\
&= & \left(\f{|v_j|^2} 4 - \f n 2 + i v_j\cdot A^\theta(x_0) \xi_j \right )u_j
  + \left( \f{|v-v_j|^2 + 2 v_j \cdot(v-v_j)} 4 \right )u_j \\
 & & + i\left( ( v-v_j)\cdot A^\theta(x_0) \xi_j\right ) u_j   + i \left( v\cdot (A^\theta(x)-A^\theta(x_0) ) \xi_j \right )u_j\\
 & &    + \left( v\cdot A^\theta(x)(\nabla_x -i \xi_j)\right)u_j\\
 &=& 0 + O(j^{\f 2 3}) + O(j^{\f 2 3}) +  O(j^{\f 3 2}) + O(j^{\f 3 2}).
\eeas
This proves
\be \label{upperb}
\| \widetilde{ P}_0(\theta) u_j\|  = O(j^{\f 3 2}).
\ee
(\ref{lowerb}) and (\ref{upperb})  contradict (\ref{nrb}) for $j$ appropriately large, because $2\sigma > \f 3 2$.
\ef

\begin{rmk}
Proposition \ref{prop.nge2} shows that for any $\ep>0$, a subelliptic estimate of the form
\[
\|H_0^{\f 3 4 +\ep} u\| +\|\w{D_x}^\ep u\| \le C_0 ( \| \widetilde{ P}_0(\theta) u\| + \| u\|), \qquad \forall u \in \vS(\bR^{2n}),
\]
is false when $\Im \theta \neq 0$ and $n \ge 2$. It also implies that $D(\widetilde{ P}_0(\theta)) \not\subseteq D(H_0)$, hence, for any $\ep>0$, $\{\widetilde{ P}_0(\theta), \theta \in \bC, |\theta| <\ep\}$ is not a holomorphic family of type (A) (\cite{k}). This shows that we cannot use the method of  analytic distortion
in the $x$-variables to define resonances of the KFP operator when $n \ge 2$.
\end{rmk}

In the remaining part of this section, we assume $n=1$.

\begin{prop} \label{prop5.2}    Let $\theta \in \bC$ with $ |\theta| < \f \pi 4$. Then $\widetilde{ P}_0(\theta)$ is  closed with domain $\vD = D(H_0)\cap D(v\cdot \nabla_x)$, and
\be \label{e5.6}
\sigma (\widetilde{ P}_0(\theta)) = \bigcup_{k\in \bN}\{k + e^{-2 i \Im \theta}r, r\ge 0\}.
\ee
\end{prop}
\pf From (\ref{e5.4}) and a density argument, one deduces that the domain of the closure of $\widetilde{ P}_0(\theta)|_{\vS}$ is equal to $\vD$.  Let $\widetilde R_0(z, \theta) = (\widetilde{ P}_0(\theta)-z)^{-1}$. For $s\ge s_0$,
\[
\widetilde R_0(-s, \theta)  -  R_0(-s, \theta) \equiv \widetilde R_0(-s, \theta) v (e^{-\theta} \nabla_x - \nabla_x^\theta)R_0(-s, \theta)
\]
modulo a term which is a compact operator. Since $n=1$,
\[
P_0(\theta) - \widetilde P_0(\theta) = v(e^{-\theta} \nabla_x - \nabla_x^\theta) + O(1) v = O(1) v\nabla_x + O(1) v,
\]
where the $O(1)$ terms are compactly supported smooth functions in $x$.
 From (\ref{sub}), one has
\[
\| v R_0(-s, \theta) \| + \| v\nabla_x  R_0(-s, \theta) \| \le C.
\]
Therefore,
$\widetilde R_0(-s, \theta) v (e^{-\theta}\na_x-\na_x^\theta)R_0(-s, \theta)\in\vL(0,0,2,1)$, hence it is a compact operator on $\vH$. The argument used in the proof of Proposition \ref{prop2.3} shows that
\[
\sigma_{\rm ess} (\widetilde P_0(\theta) ) = \sigma_{\rm ess} (P_0(\theta) )=  \bigcup_{k\in \bN}\{ k + e^{-2 i \Im \theta}r, r\ge 0\}.
\]

The proof of $ \sigma_{\rm d} (\widetilde P_0(\theta) ) = \sigma_{\rm d} (P_0(\theta) )$ (hence both are empty) is standard in the theory of resonances. We reproduce it here for the sake of completeness. The same argument is also valid to prove Proposition \ref{prop5.4} below, and will not be repeated there.   For $\theta \in \bC$ with $ |\theta | < \f \pi 4$, $\widetilde R_0(z, \theta)$ is meromorphic for $z\in \bC\setminus \sigma_{\rm ess} (\widetilde P_0(\theta) )$, which is connected.
Set
\be  \label{vA}
\vA =\{ Q(x, v)e^{-c(|x|^2 + |v|^2)}, c>0, Q \mbox{ polynomial } \}.
\ee
For $f\in \vA$, $f^\theta = U_\theta f$ and  $\widetilde f^\theta = \widetilde U_\theta f$ defined for $\theta$ real have a natural extension in $\theta\in \bC$ with $|\theta| < \f \pi 4$,  and $f^\theta$ and $\widetilde f^\theta$ belong to $\vH$. 
For $f, g\in \vA$,  $z=-s <0$ and $\theta$ real, one has
\[
\w{f, R_0(z) g} = \w{f^\theta, R_0(z, \theta)g^\theta} = \w{\widetilde f^\theta, \widetilde R_0(z, \theta)\widetilde g^\theta}.
\]
The uniqueness of meromorphic continuation shows that for $0<\Im \theta$ and $|\theta| < \f \pi 4$,
\[
\w{f^{\overline\theta}, R_0(z, \theta)g^\theta} = \w{\widetilde f^{\overline \theta}, \widetilde R_0(z, \theta)\widetilde g^\theta}
\]
as meromorphic functions of $z$ in $\bC\setminus \sigma_{\rm ess} (P_0(\theta) )$. By Proposition \ref{prop3.1}, the left-hand side of the above equality is  holomorphic for $z$ in $\bC\setminus \sigma_{\rm ess} (P_0(\theta) )$,  so is the right-hand side. Since the sets $\{\widetilde f^{\overline \theta}, f \in \vA\}$ and
$\{\widetilde g^ \theta, g\in \vA\}$ are dense in $\vH$, $\widetilde R_0(z, \theta)$ has no poles in $\bC\setminus \sigma_{\rm ess} (\widetilde P_0(\theta) )$. This proves
\[
\sigma (\widetilde{ P}_0(\theta)) = \sigma_{\rm ess} (\widetilde{ P}_0(\theta)) = \bigcup_{k\in \bN}\{k + e^{-2 i \Im \theta}r, r\ge 0\} .
\]
(\ref{e5.6}) is proved.
\ef

Assume that $V(x)$ is $C^1$ in $\bR$ and that there exist some positive constants $M, \widetilde a, C, \rho$ with $\widetilde a \le \f \pi 4$ such that $V(x)$ extends to a holomorphic function in the region
\[
\{\zeta\in \bC, |\zeta| > M, |\Im \zeta | < \tan (\tilde a) |\Re \zeta | \}
\]
and satisfies there
\be \label{ass1dst}
|V'(\zeta)| \le C \w{\zeta}^{-\rho}.
\ee
Let $\widetilde U_\theta$ be defined as before with $r_0> e^{\f \pi 4}M$. Then
\[
\widetilde P(\theta)= \widetilde U_\theta P \widetilde U_\theta^{-1}, \widetilde W_\theta = \widetilde U_\theta W\widetilde U_\theta^{-1}
\]
defined for $\theta$ real can be extended in $\theta$ into $\{\theta\in \bC, |\theta| <\widetilde a\}$. $\{\widetilde P(\theta), |\theta| < \widetilde a\}$ is a holomorphic family of type (A) in the sense of Kato, with constant domain $\vD$. Let $|\theta|<\widetilde a$ with $\Im \theta >0$. By Theorem \ref{th5.1} and Proposition \ref{prop5.2},
$ \widetilde R_0(z, \theta) \widetilde W_\theta$ is a  compact operator-valued  holomorphic function for $z\in\bC_0 \cup \Omega_-(\theta) $  and $1+ \widetilde R_0(z, \theta) \widetilde W_\theta$  is invertible for $z=-s$ with $s\ge s_0$, $s_0$ appropriately large. Similarly to Proposition \ref{prop2.3}, one has
\be
\sigma_{\rm ess}(\widetilde P(\theta)) = \bigcup_{k\in \bN}\{k + e^{-2 i \Im \theta}r, r\ge 0\}.
\ee
The analytic Fredholm theorem shows that $(1+ \widetilde R_0(z, \theta) \widetilde W_\theta)^{-1}$ is meromorphic for $z$ in $\bC_0 \cup \Omega_-(\theta)$. Making use of the equality
\[
\widetilde R(z, \theta)  =   (1+ \widetilde R_0(z, \theta) \widetilde W_\theta)^{-1}\widetilde R_0(z, \theta)
\]
which is valid for $z =-s$, $s \ge s_0$, we obtain a meromorphic continuation  of $\widetilde R(z, \theta) $ into $\bC_0 \cup \Omega_-(\theta)$.

\begin{defn} \label{defn5.3} Let $n=1$.  Assume condition (\ref{ass1dst}) for some $\widetilde a \in ]0, \f \pi 4]$. Resonances of the KFP operator $P$ defined by analytic distortion are given by
\be
\widetilde{\rm Res} (P) =\bigcup_{\theta \in \bC,\ |\theta|<\widetilde a,\ \Im \theta >0}\{ z\in \Omega_-(\theta),
\mbox{ $z$ is a pole of $\widetilde R(z, \theta)$ } \}.
\ee
\end{defn}

\begin{prop}\label{prop5.4} Assume that both conditions (\ref{ass1}) and (\ref{ass1dst}) are satisfied. Let $a_0 = \min\{a, \widetilde a\}$.  Then one has
\be
{\rm Res} (P)\cap \{z, -2 a_0 <\arg z \le 0\}= \widetilde{\rm Res} (P)\cap \{z, -2 a_0 <\arg z \le 0\}.
\ee
\end{prop}

The proof is the same as the second part of the proof of Proposition \ref{prop5.2} and is omitted.

In conclusion, this work provides a framework for further qualitative studies in physically relevant regimes, such as the large-friction and low-temperature limits. Several related questions remain open, including the definition of resonances for the KFP operator within Sj\"ostrand's theory of analytic singularities (see \cite{hs2} for resonances of elliptic operators and \cite{hss} for the tunneling effect in terms of the eigenvalues of the KFP operator).

\section*{Acknowledgments}

 Research supported in part by Nanjing University of Aeronautics and Astronautics under grant No. 90YAH25017.

\end{document}